\documentstyle{article}
\title{On the $(-1)$-curve conjecture of Friedman and Morgan.}
\author{Rogier Brussee%
\thanks{\noindent
        Math. inst. Oxford university,
        24-29 St. Giles, OX1 3LB Oxford UK, \newline
        e-mail: brussee@maths.oxford.ac.uk , fax:+44-865273583 \newline
        Supported by Nederlandse organisatie voor wetenschappelijk
        onderzoek NWO, stipend 04-63.}
}
\date{e-print alg-geom/9209001}
\catcode`\@=11
\def\[#1{\relax\ifmmode\@badmath\else\begin{equation}\label{#1}\fi}
\def\]{\relax\ifmmode\ifinner\@badmath\else\end{equation}\fi
         \else \@badmath \fi}

\def\itm#1{\par\hang\textindent{#1}\def\@currenlabel{#1}}

\newcounter{theorem}[section]  % Uses Latex to number the theorems !!!!!
\def\thetheorem{\thesection.\arabic{theorem}}

\outer\def\proclaim{\pr@claim{\bf}{ \thetheorem}{\sl}}
\outer\def\heading{\pr@claim{\bf}{ \thetheorem}{\relax}}

\outer\def\remark{\pr@claim{\sl}{ \thetheorem}{\relax}Remark }
\outer\def\acknowledgement{\pr@claim{\sl}{\relax}{\relax}Acknowledgement *.}

\def\pr@claim#1#2#3#4 #5.#6\par{\refstepcounter{theorem}\medbreak\noindent
       \def\next{#5}{#1#4#2\if\next*\else\label{\next}\fi.}\nobreak
       {\enspace#3#6}\par
  \ifdim\lastskip<\medskipamount \removelastskip\penalty55\medskip\fi}

\def\doublespace{\def\baselinestretch}

\def\Bbb{\ifmmode\let\next\Bbb@\else
 \def\next{\errmessage{Use \string\Bbb\space only in math mode}}\fi\next}
 \def\Bbb@#1{{\Bbb@@{#1}}}
 \def\Bbb@@#1{\rm\bf#1}              % thuis hack

\newif\ifcomment
\def\comment{\ifcomment\bgroup\par\medskip\noindent\small}
\def\endcomment{\par\medskip\noindent\egroup}

\newif\ifnote

\def\mathalign#1#2#3#4{\null\,\vcenter\bgroup\openup\jot\m@th
\ialign\bgroup\strut#1$\displaystyle##$#2&&#3$\displaystyle{}##$#4\crcr}
\def\endmathalign{\crcr\egroup\egroup\,}

\def\eqalign#1{\mathalign\hfil\relax\relax\hfil#1\endmathalign}
\def\cdalign#1{{\arrowlen= 4em
\mathalign\hfil\hfil\hfil\hfil#1\endmathalign}}

\let\nl=\\
\let\\=\cr

\def\eqalignno#1{\displ@y \tabskip\@centering
  \halign to\displaywidth{\hfil$\@lign\displaystyle{##}$\tabskip\z@skip
    &$\@lign\displaystyle{{}##}$\hfil\tabskip\@centering
    &\llap{$\@lign##$}\tabskip\z@skip\crcr
    #1\crcr}}

\outer\def\proof#1.{\smallbreak\noindent{\sl Proof#1\/}.\ }

\def\endproof{\nobreak{\unskip\nobreak\hfil\penalty50  % see TeX book p. 106
\hskip3em\hbox{}\nobreak\hfil$\Box$%
\parfillskip=0pt \finalhyphendemerits=0 \par
\medbreak\noindent\ignorespaces}}

\def\to{\mathchoice{\longrightarrow}{\rightarrow}%
{\rightarrow}{\rightarrow}}

\newdimen\arrowlen
\def\m@p#1#2#3#4{
       \buildrel
              \hbox spread \arrowlen{\skip@=-0.33\arrowlen plus 1 fil
                    \hskip\skip@$\m@th\scriptstyle #4$\hskip\skip@
              }
       \over{
              \mathord#1\mkern-6mu
              \cleaders\hbox{$\mkern-2mu\mathord#2\mkern-2mu$}\hfill
              \mkern-6mu\mathord#3
       }
}
\def\map{\m@p--\rightarrow}
\def\pam{\m@p\leftarrow--}
\def\equal{\m@p===}
\def\Map{\m@p==\Rightarrow}
\def\Pam{\m@p\Leftarrow==}

\def\vm@p#1#2{\Big#1
       \rlap{$\vcenter{\hbox{$\scriptstyle #2$}}$}}
\def\downmap{\vm@p\downarrow}
\def\upmap{\vm@p\uparrow}
\def\vequal{\vm@p\Vert}
\def\Downmap{\vm@p\Downarrow}

\def\inj{\lhook\nobreak\joinrel\nobreak}

\def\rmmath#1{\mathop{\rm #1}\nolimits}

\def\sheafEnd{\mathop{{\cal E}\mkern-3mu{\it nd}}\nolimits}

\def\NS{\rmmath{NS}}

\let\tensor=\otimes
\def\directsum{\mathop\oplus}

\let\slant=/

\mathcode`\:="603A %  = \colon
\mathchardef\:="303A % : as a relation

\let\iso=\cong

\def\id{{\rm id}}

\def\R{{\Bbb R}}
\def\C{{\Bbb C}}
\def\Q{{\Bbb Q}}
\def\Z{{\Bbb Z}}

\def\P{{\Bbb P}}

\def\O{{\cal O}}

\def\numfrac#1#2{\mathchoice{{\textstyle{ #1\over#2}}}%
{{ #1\over#2}}{{#1/#2}}{{#1/#2}}}
\def\quart{{\numfrac14}}
\def\half{{\numfrac12}}

\let\next=\~
\let\~=\tilde
\let\tilde=\next

\let\next=\^
\let\^=\hat
\let\hat=\next

\def\"#1{{\accent"7F \if#1i\i\else#1\fi}}

\def\({\left(}
\def\){\right)}

\def\SU{{\rm SU}}

\def\asd{{\rm asd}}
\def\Pd{{\rm Pd}}

\def\M{{\cal M}}
\def\Mbar{{\mkern4mu\overline{\mkern-4mu\M}}}

\def\<#1>{\left<#1\right>}

\catcode`\@=12

\def\Masd{\M^{\rm asd}}

\def\B{{\cal B}}
\def\diffeo{\buildrel\scriptscriptstyle \rm diff\over\iso}
\def\Pbar{{\bar\P}}
\def\closure{\overline}
\def\SO{{\rm SO}}

\def\NE{\rmmath{NE}}

\doublespace{1.1}
\rm

\fussy

\begin{document}

\maketitle

\begin{abstract}
We will prove that every differentiably embedded
sphere with self-inter\-section~$-1$ in a simply connected algebraic
surface with $p_g >0$ is homologous to an algebraic class. If the
surface has a minimal model with Picard number~$1$ or $|K_{\min}|$
contains a smooth irreducible curve of genus at least 2, and  $p_g$ is
even or $K_{\min}^2\not\equiv 7 \pmod 8$, then every such sphere is
homologous to a $(-1)$-curve, as conjectured by Friedman and Morgan.
\end{abstract}

\section{Introduction.}

It is now  well known  that the deformation type of an algebraic
surface is determined by its oriented diffeomorphism type up to a
finite number of choices \cite{F&M:ellipticI}\,\cite[theorem S.2]{F&M}.
It is therefore natural to ask if a deformation invariant is in fact an
invariant of the underlying oriented differentiable manifold. For
example,  Van de Ven conjectured that this is true for the Kodaira
dimension
\cite{O&V:overview}\,\cite{F&M}\,\cite{Pidstrigach&Tyurin:specialinst}.
In this paper we study whether the deformation invariant
decomposition
\[
{alggeodecomp} H_2(X) = H_2(X_{\min}) \oplus^\perp \directsum \Z E_i,
\]
in the homology of the minimal model and the span of the $(-1)$-curves
is invariant under orientation preserving diffeomorphisms (cf.
\cite[conj. 2]{F&M:BAMS})

A {\sl $(-1)$-curve} on a complex surface is a smooth  holomorphically
embedded 2-sphere with self-intersection $-1$. A $(-1)$-curve can be
blown down to obtain a new smooth complex surface. Successively
contracting all $(-1)$-curves gives the minimal model~$X_{\min}$ which
is unique if $p_g >0$.  More generally we will call the total transform
of a $(-1)$-curve on some intermediate blow-down a $(-1)$-curve as
well. It is in this sense that the decomposition~\ref{alggeodecomp} is
a deformation invariant.

A {\sl $(-1)$-sphere} on a $4$-manifold is  a smooth differentiably
embedded 2 sphere with self-intersection~$-1$. A classical $(-1)$-curve
is obviously a $(-1)$-sphere, a reducible one can be deformed to a
$(-1)$-sphere by smoothing out the double points.  Moreover if two
$(-1)$-curves are orthogonal, they can be deformed in disjoint
$(-1)$-spheres. Friedman and Morgan conjectured that if a surface has a
unique minimal model, then modulo homological equivalence the relation
between its $(-1)$-spheres and its $(-1)$-curves is the strongest
possible.

\proclaim Conjecture (-1)conj. (Friedman and Morgan \cite[conj. 2,3 prop.
4]{F&M:BAMS}) Let $X$ be a simply connected algebraic surface with
Kodaira dimension~$\kappa \ge 0$. Then every $(-1)$-sphere is
homologous to a $(-1)$-curve up to orientation. In particular the
decomposition~(\ref{alggeodecomp}) is invariant under orientation
preserving diffeomorphisms.

(I have slightly reformulated the conjecture, and added the simply
connectedness hypothesis). Now, $X$ contains $n$ disjoint
$(-1)$-spheres if and only if there is a differentiable connected sum
decomposition $X \diffeo X' \#n\Pbar$. The  decomposition
\ref{alggeodecomp} can be thought of as being induced by this
special connected sum decomposition.
Friedman and Morgan also made a conjecture about more general connected sum
decompositions, which would imply conjecture~\ref{(-1)conj} above.

\proclaim Conjecture consumconj. (Friedman and Morgan \cite[conj.
9]{F&M:BAMS}) Let $X$ be a simply connected algebraic surface with
$\kappa \ge 0$. Suppose $X$ admits a connected sum decomposition $X
\diffeo X' \# N$ for a negative definite manifold~$N$, then $H_2(N,\Z)$
is generated by $(-1)$-curves.

Note that by theorems of Donaldson (\cite[th. 1.3.1, 9.3.4,
10.1.1]{D&K}, $N$ has automatically a standard negative definite
intersection form if $p_g(X)>0$. Conjecture \ref{consumconj} (and hence
conjecture \ref{(-1)conj}) has been proved for blow-ups of simply
connected surfaces  with $p_g >0$ and big monodromy (like elliptic
surfaces or complete intersections), and simply connected surfaces with
$p_g >0$ whose minimal model admits a spin structure (i.e. $K_{\min}
\equiv 0(2)$) \cite[cor. 4.5.4]{F&M}. Conjecture~\ref{(-1)conj} has been
proved for the Dolgachev surfaces (i.e. $\kappa=p_g = 0$).

For minimal surfaces, conjecture \ref{(-1)conj} would imply  strong
minimality. A $4$-manifold is called {\sl strongly minimal} if for
every diffeomorphism $X \# N_1 \diffeo  Y \# N_2$ with $N_i \diffeo n_i
\Pbar^2$, we have $H_2(N_2)  \subset H_2(N_1)$ (c.f. \cite[def.
IV.4.6]{F&M}). Conjecture \ref{(-1)conj} would also imply that the canonical
class of the minimal model $K_{\min}$ is invariant mod 2 under orientation
preserving diffeomorphisms.
Conjecture \ref{consumconj} would imply that a minimal
surface is {\sl irreducible} i.e. for every decomposition $X \diffeo X'
\# N$, say $N$ is {\em homeomorphic} to  $S^4$, thereby avoiding the
Poincar\'e conjecture.
%Also, by the conjecture \ref{consumconj}, a minimal surface would be
%``very strongly minimal'' in that we may replace $N_2$ by any negative
%definite manifold.

In this paper we will show that under the stronger assumption $p_g >0$,
$(-1)$-spheres  must give rise, if not to $(-1)$-curves then in any
case to special algebraic 1-cycles. Indeed, the main theorem~\ref{main}
below leaves little room  for $(-1)$-spheres not homologous to a
$(-1)$-curve. Furthermore we will reduce a similar statement for
general connected sum decompositions, to a technical problem in gauge
theory.

To state the theorem we need some notation. Let $N_1(X)_\Z \subset
H_2(X,\Z)$ be the preferred subgroup of algebraic classes i.e. the
subgroup generated by algebraic curves. Its rank $\rho$ is  the Picard
number. The  effective cone $\NE(X) \subset N_1(X)_\Q$ is the cone
spanned by positive rational multiples of algebraic curves. The subcone
$\NE(X_{\min}) = \NE(X) \cap H_2(X_{\min},\Q)$ is the cone spanned by the
pullbacks of rational curves on the minimal model, i.e. the effective
cone in $H_2(X_{\min})$. Finally we note that since $N_1(X)_\Q$ is a
finite dimensional vector space, the closure of the effective cone
$\closure{\NE}(X)$ is well defined.

% Dually we have the Neron-Severi group $\NS(X) \subset H^2(X,\Z)$ of
%Chern classes of holomorphic line bundles, or equivalently of divisors
%mod algebraic equivalence and the nef cone, the closure (with respect
%to some norm) of the cone spanned by ample divisors in the finite
%dimensional vector space $\NS(X)_\R$. It is well known that the nef
%cone is dual to the closure of the effective cone~$\closure{\NE}(X)$
%\cite{Wilson:birational}. We can now state,

\proclaim Theorem main.  Let $X$ be  a simply connected algebraic
surface with $p_g >0$ and let~$K$ be its canonical divisor. Then for
every $(-1)$-sphere in $X$, there is an orientation such that $e$ is
either represented by a $(-1)$-curve or $e \in
\closure{\NE}(X_{\min})$, depending on whether $K\cdot e$ is  negative
or positive respectively.

Note that $K\cdot e \ne 0$ since $K\cdot e \equiv e^2 \!\!\pmod 2$.
I have no examples where  $K\cdot e$ is positive (i.e. a counter-example to
conjecture~\ref{(-1)conj}) but without further assumptions neither can I
exclude this case.

\proclaim Corollary maincor. In addition to the assumptions of the
theorem suppose that the minimal model $X_{\min}$ has Picard number $1$
or that  the linear system $|K_{\min}|$ contains a smooth irreducible
curve of genus at least two, and that $p_g$ is even or $K_{\min}^2
\not\equiv 7 \pmod 8$,
then every $(-1)$-sphere is homologous to a $(-1)$-curve (i.e.
conjecture~\ref{(-1)conj} is true for $X$).

we will use this corollary to prove conjecture~\ref{(-1)conj} for
blow-\-ups of Horikawa surfaces with $K_{\min}^2$ even and zero-sets of
general sections in  sufficiently ample $n-2$-bundles on $n$-folds with
$\rho =1$  generalising Friedman and Morgan's result for complete
intersections in $\P^n$.

The proof of theorem~\ref{main} is based on two very general properties
of the $\SO(3)$ Donaldson-Kotschick invariant~$\phi_k$. Kotschick
observed that it follows  from  the invariance properties of the
$\phi_k$~polynomial, that it is divisible by the Poincare dual of a
$(-1)$-sphere. On the other hand, using  Morgan's algebro geometric
description of the Donaldson polynomials we show that $\phi_k$  has
pure Hodge type for $k \gg 0$. The theorem and the corollary then follow
by using Donaldson's and O'Grady's non-triviality results.

\sloppy
\acknowledgement  I have greatly benefitted from discussions with
Stephan Bauer, Chris Peters, Victor Pidstrigach, Jeroen Spandaw,
Kieran O'Grady, Gang Xiao and Ping Zhang, who I would all
like to thank  heartily.  Special thanks for Simon Donaldson for
helping me with  gauge theory, and for inviting me to Oxford
university. Its  mathematics institute has proved to be a very friendly
and stimulating environment. This paper grew out of work in my thesis
\cite{RB:thesis}. It is a pleasure to thank my thesis advisors Martin
L\"ubke and  Van de Ven for their help and insight. (here  they
can not remove such words !)
\fussy

\section{The $\phi_k$~polynomials.}

We will need $\SO(3)$ Donaldson polynomials $q_{L,k,\Omega}$ and in
particular  the $\phi_k$~invariant introduced by Kotschick
\cite{Kotschick:SO(3)}. Let $X$ be a simply connected 4-manifold with
odd $b_+ \ge 3$. The polynomial $q_{L,k,\Omega}$ on $H_2(X)$
corresponds to the moduli space $\Masd(L,k)$ of ASD
$\SO(3)$-connections on the $\SO(3)$-bundle~$P_k$  with $w_2(P_k)
\equiv L \pmod 2$, and $p_1(P_k) = -4k$, oriented by the
choice of the lift~$L$ of $w_2(P_k)$, and an
orientation $\Omega$ of a maximal positive subspace in $H^2(X,\R)$
\cite[\S 9.2]{D&K}. We  choose $\Omega$ once and for all (e.g. using
a complex structure if present), and we will suppress it in the notation.
$q_{L,k}$ has degree
$
       d =  4k - \numfrac32 (1 + b_+)
$.
Note that  the $\SO(3)$ bundle~$P_k$ exists if and only if  $p_1 \equiv K^2
\pmod 4$, and that $k \in \quart \Z$.
To define $\phi_k(X)$ we lift the second
Stiefel Whitney class~$w_2(X)$ of the manifold to an integral class~$K$.
For complex surfaces, the canonical divisor is such a lift. Now define
$\phi_k = q_{K,k} $. $\phi_k(X)$ is invariant under
orientation preserving diffeomorphisms up to sign.

% Alternatively $\Masd$ can be viewed as the moduli space of projective ASD
%$U(2)$-connections on the $U(2)$ bundle with $c_1 = K$ and $c_2$
%determined by $4c_2 - K^2 = 4k$.

Now suppose that $X$  has a decomposition $X \diffeo X' \# N$ for a
negative definite  manifold $N$, necessarily with standard intersection
form. Then we have a  decomposition $H_2(X) = H_2(X') \directsum
H_2(N)$. Choose generators $e_1,\ldots,e_n$ of $H_2(N)$, such that $K\cdot
e_i \equiv -1 \pmod 4$.  This fixes the generators up to permutation.
By Poincar\'e duality we can consider the generators~$e_i$ as
linear forms on  $H_2(X)$. Any polynomial $Q$ on $H_2(X)$ can be uniquely
written as a polynomial in the dual classes of $e_i$ with polynomials
on $H_2(X')$ as coefficients ({\sl the $e$-expansion}).

\heading Definition good. A $4$-manifold is said to have a {\sl good}
connected sum decomposition $X \diffeo X' \# N$ if for every generator
$e_i$ of $H_2(N)$, $e_i$ divides $q_{L,k}(X)$ for all $L$ with
$L\cdot e_i$ odd and all $k \gg 0$.

\proclaim Proposition  {or.princip}. (Kotschick \cite[prop.
8.1]{Kotschick:SO(3)}) A connected sum decomposition $X \diffeo X' \#
n\Pbar^2$ is a good connected sum decomposition.

\proof\ (sketch). For notational simplicity we only prove divisibility
for the $\phi_k$ invariant. The reflection~$R_e$ in the hyperplane
defined by a generator~$e$ of $H_2(n\Pbar^2)$ can be represented by an
orientation preserving diffeomorphism, for example $\id \#
\C$-conjugation (c.f. \cite[prop 2.4]{F&M:ellipticI}). Then it follows
from the general invariance properties of the $\SO(3)$-polynomials
\cite[9.2.2]{D&K} that $R_e^*\phi_k(X) = - \phi_k(X)$. Hence
$\phi_k(X)$ is an odd polynomial in the dual class of $e$, in
particular $\phi_k$ is divisible by $e$.
\comment
$$
       R_e^* \phi_k = R_e^*q_{\Omega,K,k} = q_{R_e^*\Omega,R_e^*K,k}
               = (-1)^{\(K-R_e^*K \over 2\)^2}q_{\Omega,K,k} = -\phi_k.
$$
\endcomment\fi
\endproof

One should expect that any connected sum decomposition is good. This
is because the coefficients $q_{L,k,N,I}$ of the $e$-expansion have the
same invariance properties as $q_{L,k}$ under orientation preserving
diffeomorphisms of $X'$. Conjecturally, these invariants depend only on
the homotopy type of $N$, (c.f. \cite[conjecture above lemma
4.5.6]{F&M})  and so the argument for $N = n\Pbar^2$ would give the
divisibility by  the generators in general. A naive gauge theoretic
analysis seems to confirm this conjecture, but some technical
difficulties remain to be overcome. In any case we will state and prove
our results for surfaces admitting a good connected sum decomposition.

\section{Pureness of the Donaldson polynomials.}

Now we come to the algebraic geometric part of the proof. The Hodge
structure on $H^2(X,\Z)$ induces a natural Hodge structure on
$S^d H^2(X)$. Let
$$
       S^d H^2(X) \inj\map{j} H^{2d}(X\times\cdots\times X)
$$
be the natural injection in the cohomology of the $d$ fold product
of $X$. Then $j$ is a map of Hodge structures. Hence a polynomial $q\in
S^d H^2(X)$ is of pure Hodge type~$(d,d)$ if and only if $j(q)$ is pure
of Hodge type~$(d,d)$. Now clearly a sufficient condition for $j(q)$ to
be of type~$(d,d)$ is that it is represented by an algebraic cycle. We
will prove the rather natural statement that those  Donaldson
polynomials that can be  computed  completely by algebraic geometry
give rise to algebraic cycles. However to make this statement precise
requires serious work (as so often in mathematics). Fortunately almost
all of the work has already been done by J. Morgan \cite{Morgan}.

\proclaim Proposition (d,d). Let $X$ be a simply connected algebraic
surface with $p_g >0$.  Then if $L \in \NS(X)$, there is constant $k_0
>0$ such that all Donaldson polynomials $q_{L,k}$ with $k >
k_0$  and the integer  $\half (L^2 -L\cdot K) - \quart(L^2 + 4k)$ odd
are represented by algebraic cycles. In particular these polynomials
are of Hodge type $(d,d)$, where $d = \deg(q_{L,k}) = 4k - 3(1 + p_g)$.

\proof. First suppose that $L \equiv 0 \pmod 2$, then $q_{L,k}$ is up
to sign just the $\SU(2)$ polynomial $q_k$. Now the lemma follows
directly from recent results of Morgan \cite{Morgan}. He shows that for
odd $k\gg 0$, $q_k$ can be computed as follows.

Let $\Mbar^G_k = \Mbar^G_k(H,0,k)$ be the
closure of the moduli space of $H$-slope stable bundles in the moduli space
of Gieseker $H$-stable sheaves with $c_1 = 0$, $c_2 = k$. For odd~$c_2$,
there exists a universal sheaf~$\xi$ on $\Mbar^G_k$ (cf.
\cite[Remark A7]{Mukai:K3I} and \cite[prop. 2.2]{OGrady}) which
determines a correspondence
$$
\eqalign{
       \nu:H_2(X) &\to H^2(\Mbar^G_k)
\\
          \Sigma  &\to c_2(\xi) \slant \Sigma.
}
$$
Then if $H$ is $k$-generic (in a sense to be made more precise below)
and $k \gg 0$, we have \cite[theorem 1]{Morgan}
\[
{nupol}  q_k(\Sigma) = \<\nu(\Sigma)^d,[\Mbar^G_k]>.
\]

Now since the universal sheaf~$\xi$ is algebraic, it actually determines
Chow cohomology classes
$c_2(\xi) \in A^2(X\times\Mbar^G_k)$ (cf. \cite[Definition
17.3]{Fulton}). Consider the diagram
$$
\cdalign{
       X^d\times\Mbar^G_k
\\
       \llap{$\scriptstyle \pi_{X^d}$} \swarrow \quad
                      \searrow\rlap{$\scriptstyle\pi_i$}
\\
       X^d  \hskip 5em  (X\times\Mbar^G_k)_i.  \hskip -3em
}
$$
Then by equation~(\ref{nupol}), the algebraic cycle
$$
       j(q_k) = \int_{[\Mbar^G_k]} \pi_1^*c_2(\xi)\cdots \pi_d^*c_2(\xi)
\in A^d(X^d) \iso A_d(X^d),
$$
represents the image of the Donaldson polynomial~$q_k$ on the level of
Chow groups.
(Integration over the fibre $\int_{[\Mbar^G_k]}$ is
defined formally as the composition
$$
       A^i (X^d \times \Mbar^G_k) \map{[\pi_{X^d}]} A^{i-d}(X^d \times
\Mbar^G_k \to X^d) \map{\pi_*} A^{i-d}(X^d),
$$
where $[\pi_{X^d}]$  is the orientation class of the flat map
$\pi_{X^d}$ (cf. \cite[section 17.4]{Fulton}).
This proves the lemma if $L \equiv 0 \pmod 2$.

In case $L \not \equiv 0 \pmod 2$, the results of Morgan carry over
virtually unchanged, in fact the corresponding results are rather
easier. To be more precise, for a Hodge metric~$g_H$, the moduli space
of irreducible ASD $\SO(3)$-connections with $w_2 \equiv L \pmod 2$ and
$-p_1 = 4k$ can be identified with the moduli space $\M_k$ of $H$-slope
stable bundles with $c_1  = L$ and $4c_2 - c_1^2 = 4k$. For $k \gg 0$,
the closure  of the moduli space of $H$-stable bundles in the moduli
space of Gieseker stable sheaves  $\Mbar^G_k$ has the proper complex
dimension $d=4k - 3(1+p_g)$ and is generically smooth. Moreover
$\Mbar_k^G$ carries a universal sheaf~$\xi$ \cite[prop. 2.2]{OGrady},
and the class $c_2(\xi) - \quart c_1^2(\xi)$ defines a $\nu$
correspondence and Chow cohomology classes just as in the discussion
above. Finally, we choose a polarisation~$H$ which is $k$-generic in
the sense that  $H\cdot(L-2N) \ne 0$ for all $N \in \NS(X)$ with  $-4k
\le (L -2N)^2 < 0$ (i.e. $H$ is not on a wall). Then since $L \ne 0
\pmod 2$ every Gieseker $H$-semistable sheaf is actually slope
$H$-stable. Now all of the discussion in \cite{Morgan} to relate the
Gieseker and the Uhlenbeck compactification, and the $\mu$ and $\nu$
correspondence as far as it is concerned with slope stable sheaves  and
bundles carries over.
\comment
Actually the discussion of relative
classes in 6.4.2 and 6.4.3 does not carry over, but is only needed to
deal with the trivial connection anyway.
\endcomment\fi
\endproof

\remark (d,d)diff. The pureness of the polynomials
$q_{L,k}$ is also suggested by a differential geometric argument, which
seems to be the point of view taken by Tyurin
\cite[\S 4.22]{Tyurin:algaspect}. The
complex structure on $X$ induces the complex structure
$$
       T^{10}\B^*_X = \{a \in A^{10}(\sheafEnd_0(V)),\ d^* a  =0\}
$$
on the space of irreducible connections
modulo gauge $\B^*_X$. The space of irreducible ASD connections with respect
to a K\"ahler metric is then an analytic subspace. Now the explicit
formulas in \cite[Proposition 5.2.18]{D&K}
for the forms representing $\mu(\Pd(\omega))$ for an harmonic form
$\omega \in H^2(X,\C)$, show that $\mu$
preserves the Hodge structure. If we write formally
$$
       \phi_k(\Pd(\omega_1), \ldots, \Pd(\omega_d)) =
        \int_{\Masd_k}\mu(\Pd(\omega_1))\cdots \mu(\Pd(\omega_d)),
$$
then it is clear that $\phi_k(X) \ne 0$ only if the total Hodge type
of $\omega_1,\ldots,\omega_d$ is $(d,d)$. The (probably inessential)
problem is that it is not {\em a priori} obvious whether integrating
the form representatives over the non compact manifold $\M_k^\asd$
gives a valid way of computing $\phi_k$.

\section{Proof of theorem~\protect\ref{main}.}

We can now give proofs of the results stated in the introduction.
%Throughout this section we suppose that $X$ is a simply connected
%algebraic surface with $p_g >0$ admitting a good smooth connected sum
%decomposition $X \diffeo X' \# N$, for a negative definite manifold $N$
%with intersection form of rank~$n$ e.g. $X \diffeo X' \# n\Pbar^2$
The main theorem \ref{main} is the special case of theorem \ref{main'}
below for $N = n\Pbar^2$
(see definition~\ref{good} for the definition of good connected sum
decomposition).

\proclaim Theorem main'.  Let $X$ be  a simply connected algebraic
surface with $p_g >0$. Let $X_{\min}$ be its minimal model and let
$K$ be its canonical divisor.  Suppose $X$ admits a good smooth
connected sum decomposition $X \diffeo X' \# N$ for a negative definite
manifold~$N$. Then $H_2(N,\Z)$ is generated by classes
$e_1,\ldots,e_n$, with $e_i^2 = -1$ such that either $e_i$ is
represented by a $(-1)$-curve, or $e_i \in \closure{\NE}(X_{\min})$
depending on whether $K\cdot e_i$ is  negative or positive
respectively.

Here $\closure{\NE}(X_{\min})$ is the closure of the cone spanned by
positive rational multiples of algebraic curves on the minimal model

\proof. Choose a generator~$e$ of~$H_2(N)$.  We
first prove that $e$ is homologous to an algebraic cycle. Since $e$ is
certainly integral, it suffices by the Lefschetz~$(1,1)$ theorem \cite[p.
163]{G&H} to prove that its Poincar\'e dual is of
pure type~$(1,1)$.
Choose  $k$ sufficiently large as in proposition~\ref{(d,d)}
(with $L= K$), and the definition \ref{good} of good. Then $\phi_k(X)$
is non trivial and has pure Hodge type~$(d,d)$. On the other hand we
have $\phi_k = e \psi$.  Since the number of
Hodge types of~$e\psi$ is at least the number of Hodge types of~$e$,
$e$ has to be of pure type as well. Since $e$ is a real class, it is
then of type~$(1,1)$.

To show that for the proper orientation $e$ lies on the closure of the
full effective cone $\closure{\NE}(X)$,  it is enough to show that
$e\cdot H \ne 0$ for all ample divisors $H$. In fact, since the closure
of the effective cone and the nef cone are in duality \cite[proposition
2.3]{Wilson:birational}, $e \in \pm \closure{\NE}(X)$ if and only if $e$
defines a strictly positive or  strictly negative form on the ample
cone. But since the ample cone is connected, it  suffices to show that
the form~$e$ has no sign change i.e. does not vanish on the ample cone.
Since $e$ is a rational class we need to check this only for integral
ample classes. Now for a fixed ample divisor~$H$ there is a $k_0 =
k_0(H)$ such that $\phi_k(X)(H) \ne 0$ for $k > k_0$ \cite[th.
10.1.1]{D&K}. Since $e$ divides $\phi_k(X)$, it follows that $e\cdot
H\ne 0$.

By the orthogonality result \cite[th. 4.5.3]{F&M}, for every
$(-1)$-sphere~$S$ in $X$ we have either $e\cdot S = 0$ or $e = \pm
[S]$. Hence for the ``effective orientation'' of $e$ found above, $e$
is either homologous to a $(-1)$-curve, or $e$ is orthogonal to all
$(-1)$-curves, i.e. $e \in H_2(X_{\min})$. In the first case $e\cdot K
= -1 < 0$, in the latter case we have $e\cdot K = e\cdot K_{\min}> 0$
for as $p_g$ is positive, $K_{\min}$ is nef, and $K\cdot e \equiv e^2
\equiv 1 \pmod 2$. Since a divisor on the minimal model is effective if
and only if its pullback to $X$ is effective, we have
$\closure{\NE}(X_{\min})= \closure{\NE}(X) \cap H_2(X_{\min})$, and the
result follows
\endproof

Theorem \ref{main'} gives the following technical refinement of
corollary~\ref{maincor}, proving conjecture~\ref{consumconj} for
the pair $(X,N)$ under an additional hypothesis.

\proclaim Corollary maincor'. In addition to the assumptions of theorem
\ref{main'}, suppose $X$ has a deformation~$Y$ with a  minimal
model~$Y_{\min}$ such that there are no classes
$C\in\closure{\NE}(Y_{\min})$ with $C^2 = -1$ dividing all Donaldson
polynomials $q_{L,k}$ with $L\cdot C$ odd  and $k \gg 0$. Then $H^2(N)$
is generated by $(-1)$-curves. In particular this is true if
\itm{(a)}
the linear system $|K_{Y_{\min}}|$ contains a smooth irreducible curve
of genus at least 2, and either $p_g$ is even or $K_{\min}^2
\not\equiv 7 \pmod 8$, or
\itm{(b)}
the Picard number $\rho(Y_{\min}) = 1$.

As mentioned in the introduction this corollary has already been proved
without the goodness condition  under the assumptions $X_{\min}$ is spin
or $X_{\min}$ has big monodromy \cite[cor 5.4]{F&M}.

\proof. Since the deformations of a surface are all oriented
diffeomorphic, we conclude that if $X$ admits a good connected sum
decomposition $X\diffeo X' \# N$, so does its deformation~$Y$.
Moreover, the subgroup generated by $(-1)$-curves is stable under
deformation by \cite[IV.3.1]{BPV}. Hence if $H_2(N)\subset H_2(Y)$ is
generated by $(-1)$-curves, so is $H_2(N) \subset H_2(X)$. Thus we can
assume $X = Y$. By definition~\ref{good} of a good decomposition, it is
clear that no generator of $H_2(N)$ can be in
$\closure{\NE}(X_{\min})$. Hence by theorem~\ref{main'}, $H^2(N)$ is
generated by $(-1)$-curves. It remains to see that the extra condition
is satisfied in the given special cases.

In case (b), $\NS(X_{\min})$ is positive definite.  For case (a)
we argue by contradiction. Suppose there is a class
$C\in\closure{\NE}(X_{\min})$,  $C^2 = -1$ and $C$ divides $q_{C,k}$
for all $k \gg 0$.

Since $C$ is orthogonal to all $(-1)$-curves,
the proof of theorem \cite[th. 4.8]{Donaldson:pol},\,\cite[th. 9.3.14]{D&K}
gives that
$$
       q_{C,k}(X)|_{H_2(X_{\min})} = \pm q_{C,k}(X_{\min}).
$$

Since $C \in \NS(X_{\min})$, Morgan's comparison formula~\ref{nupol},
and O'Grady's non triviality result \cite[cor. 2.4, th.2.4]{OGrady},
give that for every  $\omega \in H^0(K_{\min})$, which vanishes on a smooth
irreducible curve of genus $g \ge 2$ we have
$$
       q_{C,k}(X_{\min})(\Pd(\omega + \bar \omega)) \ne 0
$$
if  $4k - 3(1+p_g)$ is even  and $k\gg 0$ with $\half(C^2-C\cdot K) -
\quart(C^2 + 4k)$ odd. (Strictly speaking O'Grady uses a slightly
different polynomial defined on $C^\perp \subset H_2(X)$, but it is
easy to see that on $C^\perp$, $q_{C,k}$ coincides with his polynomial).

Since $4k \equiv -C^2 \pmod 4$, and $\<C,\omega> =0$ this contradicts
the divisibility of $q_{C,k}$ by $C$ if $p_g$ is even. If $p_g$ is odd,
the same argument gives a contradiction  if there is a polynomial
$q_{L,k}$ with  $L\in \NS(X_{\min})$,  $L \cdot C \equiv 1$, and
$L^2\equiv K_{\min}\cdot L \equiv 0$. This an affine equation for $L
\pmod 2$ in $\NS(X_{\min}) \tensor \Z/2\Z$, so it has a solution if $C
\not\equiv K_{\min} \pmod 2$.  But if $C \equiv K_{\min}$, then
$C^2 = -1 \equiv K_{\min}^2  \pmod 8$ contrary to assumption.
\endproof

\remark *. If $|K_{\min}|$ contains a smooth
irreducible curve but $p_g$ is odd and $K_{\min}^2\equiv 7 \pmod 8$ the
proof above  shows that there is up to orientation at most one
generator $e_0$  of $H^2(N)$ which is not homologous to a $(-1)$-curve.
Hence all $(-2)$-spheres in $H_2(X_{\min})$ are orthogonal
to $e_0$, because the reflections they generate are represented by
diffeomorphisms.  We also get that if $e_0$ exists,
$w_2(X)$ is represented by the sum of the generators of $H_2(N)$, hence
$X'$ is spin.

\remark *. It follows from results in \cite{RB:e-exp} that if
$C$ is orthogonal to all $(-1)$-curves, then the divisibility of say
$\phi_k(X)$ by $C$ implies the divisibility of  $\phi_k(X')$ by $C$ for
all blow-downs of $X'$ intermediate between $X$ and $X_{\min}$. Hence
the extra condition in corollary \ref{maincor'} gets stronger as we
blow-up more points.

It is interesting to see how the minimality of a surface plays a role
here. For a non minimal surface every curve in $|K|$ is reducible since
it contains an exceptional divisor, unless $X$ is a K3-surface blown-up
once. Indeed in the non minimal case $q_{E,k}$ is divisible by the
exceptional divisor $E$, and so O'Grady's theorem could not possibly be
true. It would be very interesting if O'Grady's results would be true
assuming the existence of an irreducible curve in $|K|$, or stretching
things even further, a smooth irreducible curve in a pluricanonical
system $|nK|$.
\comment
Suppose that $X$ is of general type and
minimal.  Then say $|13K|$ contains a smooth curve of genus at least
192, and we would still have a $\theta$-characteristic that extends
over the whole surface, which seems to be one of the main points of
O'Grady's construction. Since for elliptic surfaces the $(-1)$-curve
conjecture~\ref{(-1)conj} is already known,  such a result would
confirm the conjecture for all simply connected surfaces with $p_g >0$
except those of general type with $p_g$ odd and $K_{\min}^2\equiv 7
\pmod 8$.
\endcomment\fi

\section{Examples.}

We give two examples of the use of corollary~\ref{maincor'}. The
first example follows basically by leafing through \cite{BPV}. For the
other we use Noether-Lefschetz theory  to reprove and generalise
Friedman and Morgan's result that for complete intersections with $p_g
>0$ conjecture~\ref{(-1)conj} is true. Since one approach to
Noether-Lefschetz theory is through monodromy groups, it is not
surprising that this approach gives  results similar to those using big
monodromy. However the proof may be interesting because the Noether
Lefschetz theorems we will use are proved using the ``infinitesimal
method'', based on Hodge and deformation theory.

\proclaim Proposition examples.
Suppose $X$ is a smooth simply connected surface admitting a good
connected sum decomposition $X \diffeo X' \# N$. Then $H_2(N)$ is
generated by $(-1)$-curves if
\itm(a) $X$ is the blow-up of a Horikawa surface with $K^2$
even (see \cite[table~10]{BPV}), or
\itm(b)
$(Y,\O(1))$ is a simply connected  projective local complete
intersection of dimension $r+2$ with Picard number $\rho =1$, $V$ is an
ample $r$-bundle, and $X$ is the blow up of the smooth zero locus $X_a$
of a section in $V(a)$ with $a \gg 0$.

\proof.
Case (a). By \cite[th. 10.1, remark VII.10.1]{BPV} and Bertini's
theorem, a Horikawa surface with $K^2$ even is simply connected and the
linear system $|K|$ contains a smooth curve.

\smallskip\noindent
Case (b).  For
$r=0$ the proposition is a special case of corollary~\ref{maincor'} so
we assume $r >0$.

Choose $a$ so large that $V(a)$ is globally generated.
By an application of the vector bundle version of the Lefschetz hyperplane
theorem \cite{Sommese&VdV}, \cite[cor. 22]{Okonek:Barth-Lefschetz},
$X_a$ is simply connected \cite[p.158]{G&H}.  To
prove that $p_g(X_a) >0$ consider the sequence
$$
       0 \to I_X \det(V(a))\tensor \O_Y(K_Y) \to
              \det( V(a))\tensor \O_Y(K_Y) \to \O_X(K_X) \to 0.
$$
Now choose $a$ so large that $\det V\tensor \O_Y(K_Y) \tensor \O_Y(ra)$
has a section non-vanishing on~$X$.
Finally, if $a$ is sufficiently large then $\rho(X_a) = 1$ for the
general section by Ein's generalization of the Noether-Lefschetz theorem
in $\P^3$ to ample vector bundles on projective varieties \cite[th.
2.4]{Ein}.
%Hence the proposition follows from corollary \ref{maincor'}.
\endproof

\remark *. Suppose that in case (b), $Y$
is smooth, and $V = \directsum_{i=1}^r \O(d_i)$. We can then be more
precise since there is no need to twist up. It suffices that  $V$ and
$\det(V)\tensor \O_Y(K_Y)$ are spanned by global sections,
$H^{1,1}(\det(V))= H^{1,1}(V\tensor \det(V)) = 0$,  and
$$
       H^0(V) \tensor H^0(\det(V)\tensor \O_Y(K_Y)) \to
              H^0(V\tensor \det(V)\tensor\O_Y(K_Y)))
$$
is surjective (e.g. if the $d_i \gg 0$ or $Y =\P^n$, with the
exception of $n=3$, $V = \O(2)$ or $\O(3)$, and $n=4$, $V =\O(2)
\directsum \O(2)$). This follows from judicially checking the cohomological
conditions \cite[lemma 3.2.1,3.2.2,3.2.3]{Spandaw:thesis} using
Kodaira-Nakano vanishing. For $\P^n$ the statement follows from the
classical Noether Lefschetz theorem. Also note that by choosing $a$
sufficiently large we can make $K_{X_a}$ very ample, and so we can
find a smooth irreducible curve of genus at least 2 in its linear
system.  Hence case (b) follows directly from corollary~\ref{maincor'}
if $p_g(X_a)$ or $K^2_{X_a}$ are even.

\end{document}